# Hats: all or nothing

Theo van Uem
Email: tjvanuem@gmail.com

**Abstract.**
$N$ players are randomly fitted with a colored hat ($q$ different colors). All players guess simultaneously the color of their own hat observing only the hat colors of the other $N-1$ players. The team wins if *all* players guess right. No communication of any sort is allowed, except for an initial strategy session before the game begins. In the first part of our investigation we have $q$ different colors with equal probabilities. Up to 4 colors we construct optimal strategies for any number of players using Hamming Complete Sets. For 5 colors we find optimal strategies up to 5 players using Optimal Hamming Sets. In the second part we have two colors where the probabilities may differ. We construct optimal strategies and maximal probability of winning the game for any number of players.

**Introduction.**
Hat puzzles were formulated at least since Martin Gardner's 1961 article [8]. They have got an impulse by Todd Ebert in his Ph.D. thesis in 1998 [6]. Ebert's hat problem: All players guess simultaneously the color (white or black) of their own hat observing only the hat colors of the other $N-1$ players. It is also allowed for each player to pass: no color is guessed. The team wins if at least one player guesses his hat color correctly and none of the players has an incorrect guess. Ebert's hat problem with $N=2^k-1$ players is solved in [7], using Hamming codes, and with $N=2^k$ players in [5] using extended Hamming codes. Lenstra and Seroussi [15] show that in Ebert's Hat Game, playing strategies are equivalent to binary covering codes of radius one. Ebert's asymmetric version (where the probabilities of getting a white or black hat may be different) is studied in [18],[19],[20].
In this paper $N$ distinguishable players are randomly fitted with a colored hat ($q$ different colors). All players guess simultaneously the color of their own hat observing only the hat colors of the other $N-1$ players. The team wins if *all* players guess his or her hat color correctly. An initial strategy session is allowed. Our goal is to maximize the probability of winning the game and to describe optimal strategies.
The results of our investigation are rather intriguing. For example, 10 players, 3 colors with equal probability: guessing at random gives winning probability $3^{-10}=1/59049$ , where our strategy will give probability $1/3$. Another example: 10 players, 4 colors with equal probability: guessing at random gives winning probability $4^{-10}=1/1048576$ , where our strategy will give probability $1/4$.
In case of two colors with probabilities $p$ and $q$ the maximal winning probability is: $\frac{1+|q-p|^N}{2}$ .

**PART I**

In this part $N$ distinguishable players are randomly fitted with a colored hat ($q$ different colors with equal probability).

**I.1 GOOD and BAD CASES.**
The $N$ persons in our game are distinguishable, so we can label them from 1 to $N$.
Although we are interested in the symmetric case (every color has the same probability), we will work for a great part with a generalized model: P(color $k$)= $p_k$ ($k=0,1,..,q-1$; $\sum_{k=0}^{q-1}p_k=1$).
Each possible configuration of the hats can be represented by an element of $B=\{b_1b_2\ldots b_N|b_i\in\{0,1,..,q-1\}, i=1,2..,N\}$ .
Player $i$ sees code $b_1..b_{i-1}b_{i+1}..b_N$ with decimal value $s_i=\sum_{k=1}^{i-1}b_k.q^{N-k-1}+\sum_{k=i+1}^{N}b_k.q^{N-k}$ , a value between 0 and $q^{N-1}-1$. Each player has to make a choice out of $q$ possibilities: 0='guess color 0', …, $q-1$='guess color $q-1$'.

We define a decision matrix $D = (a_{i,j})$ where $i \in \{1,2,..,N\}$ (players); $j = s_i$ ; $a_{i,j} \epsilon \{0,1,\dots,q-1\}$ (guess color 0, color 1,…, color $q-1$ ).
The meaning of $a_{i,j}$ is: player i sees j and takes decision $a_{i,j}$ (guess 0 or 1 or … $q-1$ ).
We observe the total probability (sum) of our guesses.
For each $b_1 b_2 \dots b_N$ in $B$ with $n_k$ $k$'s $(k = 0,1,..,q-1;\ \sum_{k=0}^{q-1} n_k = N)$ we have (start: sum=0):
CASE $b_1 b_2 \dots b_N$
IF $a_{1,s_1} = b_1$ AND $a_{2,s_2} = b_2$ AND … AND $a_{N,s_N} = b_N$ THEN sum=sum+$p_0^{n_0} p_1^{n_1} \dots p_{q-1}^{n_{q-1}}$ ;
(all players guess right).
We define the Hamming distance between $p_0^{n_0} p_1^{n_1} \dots p_{q-1}^{n_{q-1}}$ and $p_0^{m_0} p_1^{m_1} \dots p_{q-1}^{m_{q-1}}$ as $\sum_{i=0}^{q-1} |n_i - m_i|$, where $\sum_{i=0}^{q-1} n_i = \sum_{i=0}^{q-1} m_i = N$.
Any choice of the $a_{i,j}$ in the decision matrix determines which CASES have a positive contribution to sum (a GOOD CASE) and which CASES don't contribute positive to sum (a BAD CASE).
We focus on player $i$ $(i \in \{1,..,N\})$. Each $a_{i,j} = m$ $(m \in \{0,1,\dots,q-1\})$ has $(q-1)$ counterparts $a_{i,j} = k$ $(k \in \{0,1,\dots,q-1\} \backslash \{m\})$: use the flipping procedure in position $i$:
CASE $b_1 .. b_{i-1} m b_{i+1} .. b_N \rightarrow$ CASE $b_1 .. b_{i-1} k b_{i+1} .. b_N$. So, when the GOOD CASE has probability $p_0^{n_0} p_1^{n_1} \dots p_{q-1}^{n_{q-1}}$ we get, for fixed $m$ and $k$, a BAD CASE with probability $p_1^{n_1} p_2^{n_2} \dots p_q^{n_q} p_m^{-1} p_k$ and Hamming distance 2 to our GOOD CASE.
We can't have two GOOD CASES G1 and G2 with Hamming distance 2: if the Hamming distance is 2, then we can get G2 by a bit flip of G1, but this changes the process $a_{1,s_1} = b_1$ AND $a_{2,s_2} = b_2$ AND … AND $a_{N,s_N} = b_N$ and we get a BAD CASE G2.
In the next sections we renumber the colors 1,2,…,q instead of 0,1,…,q-1.

### I.2 Hamming Complete Set (HCS)

$\mathcal{H} = \{S_1, S_2, .., S_q\}$ where $S_i \subset \sum_{\sum_{m=1}^{q} n_m = N} (p_1^{n_1} \dots p_q^{n_q})$ $(i = 1,2,..,q)$.
$\mathcal{H}$ is Hamming Complete when we have the following three properties:

- (i) $S_i$ consists of elements with Hamming distance greater than 2
- (ii) The sets $S_i$ are disjunct
- (iii) Completeness: $\bigcup_i S_i = \sum_{\sum_{m=1}^{q} n_m = N} (p_1^{n_1} \dots p_q^{n_q})$

If we can construct a Hamming Complete Set (HCS), then each $S_i$ induces a strategy $\mathcal{S}_i$ with probability $\mathcal{P}_i$ (this will become clear in the next sections).
We now return to the symmetric case: each color has probability $1/q$. This is an upper bound of $\mathcal{P}_i$ (we can't do better than the result of one player), so we have: $\mathcal{P}_i \leq 1/q$ $(i = 1,2,..,q)$ and $\sum_{i=1}^{q} \mathcal{P}_i = 1$.
***Conclusion: All strategies $\boldsymbol{\mathcal{S}_i}$ $(\boldsymbol{i = 1,2,..,q})$ are optimal and the probability of each strategy is $\boldsymbol{1/q}$, independent of the number of players.***
In the next sections we construct Hamming Complete Sets up to 4 colors.

### I.3 Two color Hat Game

Hamming Complete Set:
$S_1 = \sum_{k\ even} p_1^{N-k} {p_2}^k$
$S_2 = \sum_{k\ odd} p_1^{N-k} {p_2}^k$
$\boldsymbol{\mathcal{S}_1}$: guess in such a way that there is an even number of hats of color 2.
$\boldsymbol{\mathcal{S}_2}$: guess in such a way that there is an odd number of hats of color 2.
(All players have to make the same choice in the initial strategy session before the game begins).
$\mathcal{P}_1 = \mathcal{P}_2 = 1/q = 1/2$ , which agrees with
$\mathcal{P}_1 = \sum_{k\ even} \binom{N}{k} {p_1}^{N-k} {p_2}^k$ and $\mathcal{P}_2 = \sum_{k\ odd} \binom{N}{k} {p_1}^{N-k} {p_2}^k$ where $p_1 = p_2 = 1/2$.

## I.4 Three color Hat Game

Operator $T$ is defined by $T(p_2^m p_3^n) = \sum_{s=-[\frac{n}{3}]}^{[\frac{m}{3}]} p_2^{m-3s} p_3^{n+3s}$.

$[p_1 + (p_2 + p_3)]^N = \sum_{k=0}^{N} \binom{N}{k} p_1{}^{N-k} \, (p_2 + p_3)^k = 1$ .

Concentrating on $(p_2 + p_3)^k$, we construct a Hamming Complete Set:

$$S_1 = \sum_{n=0}^{[\frac{N}{2}]} p_1^{N-2n} T(p_2^n p_3^n) + \sum_{n=0}^{[\frac{N-3}{2}]} p_1^{N-2n-3} T(p_2^{n+3} p_3^n)$$

$$S_2 = \sum_{n=0}^{[\frac{N-1}{2}]} p_1^{N-2n-1} T(p_2^{n+1} p_3^n) + \sum_{n=0}^{[\frac{N-2}{2}]} p_1^{N-2n-2} T(p_2^n p_3^{n+2})$$

$$S_3 = \sum_{n=0}^{[\frac{N-1}{2}]} p_1^{N-2n-1} T(p_2^n p_3^{n+1}) + \sum_{n=0}^{[\frac{N-2}{2}]} p_1^{N-2n-2} T(p_2^{n+2} p_3^n)$$

***Conclusion: All strategies $\mathcal{S}_i$ $(i = 1, 2, .., 3)$ are optimal and the probability of each strategy is $1/3$, independent of the number of players.***

We notice that $\mathcal{S}_3$ can be found by interchanging in $\mathcal{S}_2$ the colors 2 and 3.

We illustrate our theory with three examples.

*Example 1*: Three players and three colors.
$S_1 = \sum_{n=0}^{[\frac{N}{2}]} p_1^{N-2n} T(p_2^n p_3^n) + \sum_{n=0}^{[\frac{N-3}{2}]} p_1^{N-2n-3} T(p_2^{n+3} p_3^n) = p_1^3 T(1) + p_1 T(p_2 p_3) + T(p_2^3) = p_1^3 + p_1 p_2 p_3 + (p_2^3 + p_3^3)$ inducing strategy $\mathcal{S}_1$ = {300, 111, 030, 003}: when you see two identic colors, guess that color, otherwise choose the missing color.
The probability is $p_1^3 + 6p_1p_2p_3 + (p_2^3 + p_3^3) = 1/3$.
$S_2 = \sum_{n=0}^{[\frac{N-1}{2}]} p_1^{N-2n-1} T(p_2^{n+1} p_3^n) + \sum_{n=0}^{[\frac{N-2}{2}]} p_1^{N-2n-2} T(p_2^n p_3^{n+2}) = p_1^2 T(p_2) + T(p_2^2 p_3) + p_1 T(p_3^2) = p_1^2 p_2 + p_2^2 p_3 + p_1 p_3^2$ with probability $3p_1^2 p_2 + 3p_2^2 p_3 + 3p_1 p_3^2$, inducing strategy $\mathcal{S}_2$ ={210, 021, 102}.
$\mathcal{S}_3 = \{201, 012, 120\}$.

*Example 2*: Four players and three colors.
$S_1 = \sum_{n=0}^{[\frac{N}{2}]} p_1^{N-2n} T(p_2^n p_3^n) + \sum_{n=0}^{[\frac{N-3}{2}]} p_1^{N-2n-3} T(p_2^{n+3} p_3^n) = p_1^4 T(1) + p_1^2 T(p_2 p_3) + T(p_2^2 p_3^2) + p_1 T(p_2^3) = p_1^4 + p_1^2 p_2 p_3 + p_2^2 p_3^2 + p_1 (p_2^3 + p_3^3)$ , probability $p_1^4 + 12 p_1^2 p_2 p_3 + 6 p_2^2 p_3^2 + 4 p_1 (p_2^3 + p_3^3) = 1/3$ ; strategy $\mathcal{S}_1$ ={400, 211, 022, (130, 103)}.
$S_2 = \sum_{n=0}^{[\frac{N-1}{2}]} p_1^{N-2n-1} T(p_2^{n+1} p_3^n) + \sum_{n=0}^{[\frac{N-2}{2}]} p_1^{N-2n-2} T(p_2^n p_3^{n+2}) = p_1^3 T(p_2) + p_1 T(p_2^2 p_3) + p_1^2 T(p_3^2) + T(p_2 p_3^3) = p_1^3 p_2 + p_1 p_2^2 p_3 + p_1^2 p_3^2 + (p_2 p_3^3 + p_2^4)$, probability $4p_1^3 p_2 + 12 p_1 p_2^2 p_3 + 6 p_1^2 p_3^2 + 4 p_2 p_3^3 + p_2^4$;
$\mathcal{S}_2 = \{320, 121, 202, 013, 040\}$ ; $\mathcal{S}_3 = \{302, 112, 220, 031, 004\}$

*Example 3*: Ten players and three colors.
We limit to:
$S_1 = \sum_{n=0}^{[\frac{N}{2}]} p_1^{N-2n} T(p_2^n p_3^n) + \sum_{n=0}^{[\frac{N-3}{2}]} p_1^{N-2n-3} T(p_2^{n+3} p_3^n) = p_1^{10} T(1) + p_1^8 T(p_2 p_3) + p_1^6 T(p_2^2 p_3^2) + p_1^4 T(p_2^3 p_3^3) + p_1^2 T(p_2^4 p_3^4) + T\left(p_2^5 p_3^5\right) + p_1^7 T(p_2^3) + p_1^5 T(p_2^4 p_3) + p_1^3 T\left(p_2^5 p_3^2\right) + p_1 T(p_2^6 p_3^3) = p_1^{10} + p_1^8 p_2 p_3 + p_1^6 p_2^2 p_3^2 + p_1^4 (p_2^6 + p_2^3 p_3^3 + p_3^6) + p_1^2 (p_2^7 p_3 + p_2^4 p_3^4 + p_2 p_3^7) + \left(p_2^8 p_3^2 + p_2^5 p_3^5 + p_2^2 p_3^8\right) + p_1^7 (p_2^3 + p_3^3) + p_1^5 \, (p_2^4 p_3 + p_2 p_3^4) + p_1^3 \, ( p_2^5 p_3^2 + p_2^2 p_3^5) + p_1 (p_2^9 + p_2^6 p_3^3 + p_2^3 p_3^6 + p_3^9)$
Strategy $\mathcal{S}_1$ ={10.0.0, 811, 622, (460, 433, 406), (271, 244, 217), (082, 055, 028), (730, 703), (541, 514), (352, 325), (190, 163, 136, 109)}.

There is a way to determine the optimal strategies without using the operator $T$.
We give the key for strategy $\mathcal{S}_1$:

$$\sum_{k\ even}\binom{N}{k}\sum_{s=0}^{[\frac{k}{3}]}\binom{k}{3s+k_1}+\sum_{k\ odd}\binom{N}{k}\sum_{s=0}^{[\frac{k}{3}]}\binom{k}{3s+k_2}=3^{N-1}$$

where $k_1=\frac{k}{2}mod3$ and $k_2=\frac{(k-3)}{2}mod3$

Taking $N$ =3, we get: $\binom{3}{0}\binom{0}{0}+\binom{3}{2}\binom{2}{1}+\binom{3}{3}\{\binom{3}{0}+\binom{3}{3}\}$, inducing strategy {300, 111, 030, 003}.

When $N$ =4 we get: $\binom{4}{0}\binom{0}{0}+\binom{4}{2}\binom{2}{1}+\binom{4}{4}\binom{4}{2}+\binom{4}{3}\{\binom{3}{0}+\binom{3}{3}\}$, inducing strategy {400, 211, 022, 130, 103}.

When $N=10$ we get:

$$\binom{10}{0}\binom{0}{0}+\binom{10}{2}\binom{2}{1}+\binom{10}{4}\binom{4}{2}+\binom{10}{6}\left\{\binom{6}{0}+\binom{6}{3}+\binom{6}{6}\right\}+..+\binom{10}{9}\left\{\binom{9}{0}+\binom{9}{3}+\binom{9}{6}+\binom{9}{9}\right\}.$$

$\mathcal{S}_1=$ {10.0.0, 811, 622, (460, 433, 406), …, ( 190, 163, 136, 109)}.

We can obtain similar results, without using the operator $T$, in case of $\mathcal{S}_2$ and $\mathcal{S}_3$:
We give the key for strategy $\mathcal{S}_2$:

$$\sum_{k\ even}\binom{N}{k}\sum_{s=0}^{[\frac{k}{3}]}\binom{k}{3s+k_3}+\sum_{k\ odd}\binom{N}{k}\sum_{s=0}^{[\frac{k}{3}]}\binom{k}{3s+k_4}=3^{N-1}$$

where $k_3=\frac{(k+2)}{2}mod3$ and $k_4=\frac{(k-1)}{2}mod3$

The key for strategy $\mathcal{S}_3$:

$$\sum_{k\ even}\binom{N}{k}\sum_{s=0}^{[\frac{k}{3}]}\binom{k}{3s+k_5}+\sum_{k\ odd}\binom{N}{k}\sum_{s=0}^{[\frac{k}{3}]}\binom{k}{3s+k_6}=3^{N-1}$$

where $k_5=\frac{(k-2)}{2}mod3$ and $k_6=\frac{(k+1)}{2}mod3$

## I.5 Four color Hat Game

$$[p_1+(p_2+p_3+p_4)]^N=\sum_{k=0}^{N}\binom{N}{k}{p_1}^{N-k}\,(p_2+p_3+p_4)^k=1$$

Concentrating on $(p_2+p_3+p_4)^k$, we construct a Hamming Complete Set:

$$S_1=\sum_{n=0}^{[\frac{N}{2}]}p_1^{N-2n}(p_2^2+p_3^2+p_4^2)^{n*}+\sum_{n=0}^{[\frac{N-3}{2}]}p_1^{N-2n-3}p_2p_3p_4(p_2^2+p_3^2+p_4^2)^{n*}$$

$$S_2=\sum_{n=0}^{[\frac{N-1}{2}]}p_1^{N-2n-1}p_2(p_2^2+p_3^2+p_4^2)^{n*}+\sum_{n=0}^{[\frac{N-2}{2}]}p_1^{N-2n-2}p_3p_4(p_2^2+p_3^2+p_4^2)^{n*}$$

$$S_3=\sum_{n=0}^{\left[\frac{N-1}{2}\right]}p_1^{N-2n-1}p_3(p_2^2+p_3^2+p_4^2)^{n*}+\sum_{n=0}^{\left[\frac{N-2}{2}\right]}p_1^{N-2n-2}p_2p_4(p_2^2+p_3^2+p_4^2)^{n*}$$

$$S_4=\sum_{n=0}^{[\frac{N-1}{2}]}p_1^{N-2n-1}p_4(p_2^2+p_3^2+p_4^2)^{n*}+\sum_{n=0}^{[\frac{N-2}{2}]}p_1^{N-2n-2}p_2p_3(p_2^2+p_3^2+p_4^2)^{n*}$$

where we define $(p_2^2+p_3^2+p_4^2)^{n*}=\sum_{i+j+k=n}p_2^{2i}\,p_3^{2j}p_4^{2k}$.

(i), (ii) are easily verified.
To prove completeness, we first consider the sum of coefficients of $(p_2+p_3+p_4)^{2k}$:

$$\sum_{i+j+l=k}\binom{2k}{2i,2j,2l}+3\sum_{i+j+l=k-1}\binom{2k}{2i,2j+1,2l+1}=(1+1+1)^{2k}=3^{2k}$$

The sum of coefficients of $(p_2+p_3+p_4)^{2k+1}$:

$$\sum_{i+j+l=k-1}\binom{2k+1}{2i+1,2j+1,2l+1}+3\sum_{i+j+l=k}\binom{2k+1}{2i,2j,2l+1}=3^{2k+1}$$

***Conclusion: All strategies $S_i$ $(i = 1, 2, .., 4)$ are optimal and the probability of each strategy is $1/4$, independent of the number of players.***

We notice that $\mathcal{S}_j$ can be found by interchanging in $\mathcal{S}_k$ the colors $j$ and $k$ $(j \neq k, j \geq 2, k \geq 2)$.

We give some examples.

*Example 1*: Two players, four colors.

$S_1 = \sum_{n=0}^{[\frac{N}{2}]} p_1^{N-2n}(p_2^2 + p_3^2 + p_4^2)^{n*} + \sum_{n=0}^{[\frac{N-3}{2}]} p_1^{N-2n-3} p_2p_3p_4(p_2^2 + p_3^2 + p_4^2)^{n*} = p_1^2 + p_2^2 + p_3^2 + p_4^2$

, which is also the probability (1/4), inducing strategy: {2000,0200,0020,0002}: guess what you see.

$S_2 = \sum_{n=0}^{[\frac{N-1}{2}]} p_1^{N-2n-1} p_2(p_2^2 + p_3^2 + p_4^2)^{n*} + \sum_{n=0}^{[\frac{N-2}{2}]} p_1^{N-2n-2} p_3p_4(p_2^2 + p_3^2 + p_4^2)^{n*} = p_1p_2 + p_3p_4;$

probability $2p_1p_2 + 2p_3p_4$; strategy: {1100, 0011}.

$\mathcal{S}_3 =$ {1010, 0101}; $\mathcal{S}_4 =$ {1001, 0110}

In the next examples we only consider $S_1$.

*Example 2*: Three players, four colors.

$S_1 = \sum_{n=0}^{[\frac{N}{2}]} p_1^{N-2n}(p_2^2 + p_3^2 + p_4^2)^{n*} + \sum_{n=0}^{[\frac{N-3}{2}]} p_1^{N-2n-3} p_2p_3p_4(p_2^2 + p_3^2 + p_4^2)^{n*} = p_1^3 + p_1(p_2^2 + p_3^2 + p_4^2) + p_2p_3p_4$ ; probability $p_1^3 + 3p_1(p_2^2 + p_3^2 + p_4^2) + 6p_2p_3p_4$ ;

$\mathcal{S}_1$ ={3000, (1200, 1020, 1002), 0111}.

*Example 3*: Four players, four colors.

$S_1 = \sum_{n=0}^{[\frac{N}{2}]} p_1^{N-2n}(p_2^2 + p_3^2 + p_4^2)^{n*} + \sum_{n=0}^{[\frac{N-3}{2}]} p_1^{N-2n-3} p_2p_3p_4(p_2^2 + p_3^2 + p_4^2)^{n*} = p_1^4 + p_1^2(p_2^2 + p_3^2 + p_4^2) + (p_2^2 + p_3^2 + p_4^2)^{2*} + p_1p_2p_3p_4$ ;

probability $p_1^4 + 6p_1^2(p_2^2 + p_3^2 + p_4^2) + (p_2^4 + p_3^4 + p_4^4) + 6(p_2^2p_3^2 + p_4^2p_2^2 + p_3^2p_4^2) + 24p_1p_2p_3p_4$.

$\mathcal{S}_1 =$ {4000, (2200,2020, 2002),(0400, 0040, 0004,0220,0202,0022), 1111}

*Example 4*: Five players, four colors.

$S_1 = \sum_{n=0}^{[\frac{N}{2}]} p_1^{N-2n}(p_2^2 + p_3^2 + p_4^2)^{n*} + \sum_{n=0}^{[\frac{N-3}{2}]} p_1^{N-2n-3} p_2p_3p_4(p_2^2 + p_3^2 + p_4^2)^{n*} = p_1^5 + p_1^3(p_2^2 + p_3^2 + p_4^2) + p_1\ (p_2^2 + p_3^2 + p_4^2)^{2*} + p_1^2p_2p_3p_4 + p_2p_3p_4(p_2^2 + p_3^2 + p_4^2)$ .

$\mathcal{S}_1$ ={ 5000, (3200, 3020, 3002), (1400, 1040, 1004, 1220, 1202, 1022), 2111, (0311, 0131, 0113)}.

***When we have more than 4 colors, then we can't construct HCS and we need a new concept: Optimal Hamming Set.***

## I.6 Optimal Hamming Set

We have N players and each player is randomly fitted with 1 of q different colors. Let $n_i$ be the number of players with color $i$ $(i = 1,2, .., q)$. We describe this situation by $[n_1, n_2, \dots, n_q]$ , where $\sum_{i=1}^{q} n_i = N$.

The Hamming distance between $[n_1, n_2, \dots, n_q]$ and $[m_1, m_2, \dots, m_q]$ is defined by $\sum_{i=1}^{q} |\, n_i - m_i|$.

An Optimal Hamming Set (OHS) consists of elements $[n_1, n_2, \dots, n_q]$ with two restrictions:

(i) All elements have Hamming distance greater than two to each other

(ii) $\sum_{OHS} \frac{N!}{n_1!n_2!\dots n_q!} = q^{N-1}$

The motivation for restriction (i) can be found in section I.1.

Restriction (ii) implies a probability $\frac{q^{N-1}}{q^N} = \frac{1}{q}$ in the outcome space where each of the N players has a choice out of $q$ different colors. And the probability of winning our game can't exceed $\frac{1}{q}$.

In Appendix A we give results of optimal covering with 5 colors up to 5 players.

**PART II**

In this part $N$ distinguishable players are fitted at random with a white or black hat, where the probabilities of getting a white or black hat ($p$ respectively $q$; $p + q = 1$ ) may be different, but known and the same to all players.

**Lemma II.1**
The sequence $\binom{N}{k} p^{N-k} q^k$ $(k = 0,1,..,N)$ consists of two parts: the first part is monotone increasing (may be empty) and the second part is monotone decreasing (may be empty).
*Proof*
Increasing values for $k = 1,2,..,s$ will be found if: $\frac{\binom{N}{k} p^{N-k} q^k}{\binom{N}{k-1} p^{N-k+1} q^{k-1}} = \frac{(N-k+1)}{k}\frac{q}{p} > 1$, so $p < \frac{N-(k-1)}{N+1}$. It follows: $p < \frac{N-(s-1)}{N+1}$. In the same way we can show that decreasing values for $k = s+1, s+2,..,N$ comes to $p > \frac{N-(k-1)}{N+1}$. So: $p > \frac{N-s}{N+1}$.
$\forall_{p\in(0,1)} \exists_{s\in\{0,1,..,N\}}$: $\frac{N-s}{N+1} < p < \frac{N-(s-1)}{N+1}$ and for this value of $s$ we have:
sequence is increasing if $k = 1,2,..,s$ and decreasing if $k = s+1, s+2,..,N$. ∎

We use a result of I.1: GOOD CASES have Hamming distance greater than 2 to each other.
An optimal solution (maximal probability) can be written as:
$\sum_{k=k_1,k_2,..,k_r} \binom{N}{k} p^{N-k} q^k$ for some $r$ and $0 \le k_1 < k_2 <..< k_r \le N$.
We define the *gap* between $k_{i+1}$ and $k_i$ as $k_{i+1} - k_i$ $(i = 1,..,r-1)$.
A *gap g* between $k_{i+1}$ and $k_i$ correspondents with Hamming distance *2g* between $p^{N-k_i} q^{k_i}$ and $p^{N-k_{i+1}} q^{k_{i+1}}$
Hamming distance is greater than 2, so: $k_{i+1} - k_i > 1$ $(i = 1,..,r-1)$.
The maximal *gap* in an optimal solution can't be 3 or more. If the *gap* is 3 or more, then we can get a better solution by the following receipt:
If $k_i$ and $k_{i+1}$ are in the monotone increasing part of $\binom{N}{k} p^{N-k} q^k$, then replace $k_i$ by $k_i + 1$.
If $k_i$ and $k_{i+1}$ are in the monotone decreasing part of $\binom{N}{k} p^{N-k} q^k$, then replace $k_{i+1}$ by $k_{i+1} - 1$.
If $k_i$ is in the monotone increasing part and $k_{i+1}$ is in the monotone decreasing part of $\binom{N}{k} p^{N-k} q^k$, then we can always push up at least one $k$ $(k_i \le k \le k_{i+1})$.
So, we have:

**Lemma II.2**
The *gap* in an optimal solution is always 2. ∎

We define 0-parity of a CASE as the parity of the number of white hats in that CASE.

**Theorem II.1**
The maximal probabilities and the optimal strategies of our hat problem are:
when $p < q$ or $N$ is even:

$\sum_{k\ even} \binom{N}{k} p^k q^{N-k} = \frac{1+(q-p)^N}{2}$; strategy: even 0-parity,

when $p > q$ and $N$ is odd:

$\sum_{k\ odd} \binom{N}{k} p^k q^{N-k} = \frac{1-(q-p)^N}{2}$; strategy: odd 0-parity,

when $p = q$:

$\frac{1}{2}$; strategy: all players odd 0-parity or all players even 0-parity.

*Proof*

Lemma II.2 gives that an optimal solution has all even or all odd terms in $\sum_k \binom{N}{k} p^k q^{N-k}$.

We have: $\sum_{k\ even} \binom{N}{k} p^k q^{N-k} - \sum_{k\ odd} \binom{N}{k} p^k q^{N-k} = (q-p)^N$.

Using $\sum_{k\ even} \binom{N}{k} p^k q^{N-k} + \sum_{k\ odd} \binom{N}{k} p^k q^{N-k} = (q+p)^N = 1$, we obtain our goal. ■

So, we have: *maximal winning probability:* $\frac{1+|q-p|^N}{2}$ $(N = 1,2,3,\dots)$.

***Hats: Nothing.***

Until now we were interested in perfect guessing: all players must guess correct. In this section we demand the opposite: *all players must guess wrong*.

The theory of maximal winning probabilities stays the same, but the optimal strategy is now: use the strategy in 'all players guess right', followed by a bit flip.

**Appendix A (Part I) Optimal covering with 5 colors up to 5 players.**

**B.1 Five colors and two players.**

The optimal strategy is simple: guess the same color as you see.

In the next table we see OHS, $\sum_{OHS} \frac{N!}{n_1!n_2!\dots n_q!} = q^{N-1}$ and Hamming distances:

| | | | | | |
|---|---|---|---|---|---|
| 2 | 0 | 0 | 0 | 0 | 1 |
| 0 | 2 | 0 | 0 | 0 | 1 |
| 0 | 0 | 2 | 0 | 0 | 1 |
| 0 | 0 | 0 | 2 | 0 | 1 |
| 0 | 0 | 0 | 0 | 2 | 1 |
| | | | | | 5 |

| | 1 | 2 | 3 | 4 |
|---|---|---|---|---|
| 1 | | | | |
| 2 | 4 | | | |
| 3 | 4 | 4 | | |
| 4 | 4 | 4 | 4 | |
| 5 | 4 | 4 | 4 | 4 |

**B.2 Five colors and three players.**

In the next table we see OHS, $\sum_{OHS} \frac{N!}{n_1!n_2!\dots n_q!} = q^{N-1}$ and Hamming distances:

| | | | | | |
|---|---|---|---|---|---|
| 3 | 0 | 0 | 0 | 0 | 1 |
| 1 | 1 | 1 | 0 | 0 | 6 |
| 1 | 0 | 0 | 1 | 1 | 6 |
| 0 | 2 | 0 | 1 | 0 | 3 |
| 0 | 0 | 2 | 0 | 1 | 3 |
| 0 | 0 | 1 | 2 | 0 | 3 |
| 0 | 1 | 0 | 0 | 2 | 3 |
| | | | | | 25 |

| | 1 | 2 | 3 | 4 | 5 | 6 |
|---|---|---|---|---|---|---|
| 1 | | | | | | |
| 2 | 4 | | | | | |
| 3 | 4 | 4 | | | | |
| 4 | 6 | 4 | 4 | | | |
| 5 | 6 | 4 | 4 | 6 | | |
| 6 | 6 | 4 | 4 | 4 | 4 | |
| 7 | 6 | 4 | 4 | 4 | 4 | 6 |

The strategy is: observe the other two players; when you see color 1 two times then guess color 1 (first row); when you see color 2 two times then guess color 4 (row four); etc.

By interchanging two colors (two rows) we get another OHS, but this new solution is not disjunct with our first solution (disjunction was demanded by Hamming Complete Set, which works only up to 4 colors).

**B.3 Five colors and four players.**

In the next tables we see OHS, $\sum_{OHS} \frac{N!}{n_1!n_2!\dots n_q!} = q^{N-1}$ and Hamming distances:

| | | | | | | |
|---|---|---|---|---|---|---|
| 1 | 4 | 0 | 0 | 0 | 0 | 1 |
| 2 | 0 | 4 | 0 | 0 | 0 | 1 |
| 3 | 0 | 0 | 4 | 0 | 0 | 1 |
| 4 | 0 | 0 | 0 | 4 | 0 | 1 |
| 5 | 0 | 0 | 0 | 0 | 4 | 1 |
| 6 | 2 | 1 | 1 | 0 | 0 | 12 |
| 7 | 2 | 0 | 0 | 1 | 1 | 12 |
| 8 | 1 | 2 | 0 | 1 | 0 | 12 |
| 9 | 0 | 2 | 1 | 0 | 1 | 12 |
| 10 | 1 | 0 | 2 | 0 | 1 | 12 |
| 11 | 0 | 1 | 2 | 1 | 0 | 12 |
| 12 | 1 | 0 | 1 | 2 | 0 | 12 |
| 13 | 0 | 1 | 0 | 2 | 1 | 12 |
| 14 | 1 | 1 | 0 | 0 | 2 | 12 |
| 15 | 0 | 0 | 1 | 1 | 2 | 12 |
| | | | | | | 125 |

| | 1 | 2 | 3 | 4 | 5 | 6 | 7 | 8 | 9 | 10 | 11 | 12 | 13 | 14 |
|---|---|---|---|---|---|---|---|---|---|---|---|---|---|---|
| 2 | 8 | | | | | | | | | | | | | |
| 3 | 8 | 8 | | | | | | | | | | | | |
| 4 | 8 | 8 | 8 | | | | | | | | | | | |
| 5 | 8 | 8 | 8 | 8 | | | | | | | | | | |
| 6 | 4 | 6 | 6 | 8 | 8 | | | | | | | | | |
| 7 | 4 | 8 | 8 | 6 | 6 | 4 | | | | | | | | |
| 8 | 6 | 4 | 8 | 6 | 8 | 4 | 4 | | | | | | | |
| 9 | 8 | 4 | 6 | 8 | 6 | 4 | 6 | 4 | | | | | | |
| 10 | 6 | 8 | 4 | 8 | 6 | 4 | 4 | 6 | 4 | | | | | |
| 11 | 8 | 6 | 4 | 6 | 8 | 4 | 6 | 4 | 4 | 4 | | | | |
| 12 | 6 | 8 | 6 | 4 | 8 | 4 | 4 | 4 | 6 | 4 | 4 | | | |
| 13 | 8 | 6 | 8 | 4 | 6 | 6 | 4 | 4 | 4 | 6 | 4 | 4 | | |
| 14 | 6 | 6 | 8 | 8 | 4 | 4 | 4 | 4 | 4 | 4 | 6 | 6 | 4 | |
| 15 | 8 | 8 | 6 | 6 | 4 | 6 | 4 | 6 | 4 | 4 | 4 | 4 | 4 | 4 |

### B.4 Five colors and five players.

In the next table we see OHS, $\sum_{OHS} \frac{N!}{n_1!n_2!\ldots n_q!} = q^{N-1}$ and Hamming distances:

| | | | | | | |
|---|---|---|---|---|---|---|
| 1 | 2 | 1 | 1 | 1 | 0 | 60 |
| 2 | 1 | 2 | 1 | 0 | 1 | 60 |
| 3 | 1 | 0 | 2 | 1 | 1 | 60 |
| 4 | 0 | 1 | 1 | 2 | 1 | 60 |
| 5 | 1 | 1 | 0 | 1 | 2 | 60 |
| 6 | 2 | 0 | 0 | 2 | 1 | 30 |
| 7 | 2 | 0 | 1 | 0 | 2 | 30 |
| 8 | 0 | 2 | 2 | 1 | 0 | 30 |
| 9 | 1 | 2 | 0 | 2 | 0 | 30 |
| 10 | 0 | 1 | 2 | 0 | 2 | 30 |
| 11 | 3 | 1 | 0 | 0 | 1 | 20 |
| 12 | 0 | 3 | 0 | 1 | 1 | 20 |
| 13 | 1 | 1 | 3 | 0 | 0 | 20 |
| 14 | 1 | 0 | 1 | 3 | 0 | 20 |
| 15 | 0 | 0 | 1 | 1 | 3 | 20 |
| 16 | 3 | 0 | 2 | 0 | 0 | 10 |
| 17 | 2 | 3 | 0 | 0 | 0 | 10 |
| 18 | 0 | 0 | 3 | 2 | 0 | 10 |
| 19 | 0 | 0 | 0 | 3 | 2 | 10 |
| 20 | 0 | 2 | 0 | 0 | 3 | 10 |
| 21 | 4 | 0 | 0 | 1 | 0 | 5 |
| 22 | 0 | 4 | 1 | 0 | 0 | 5 |
| 23 | 0 | 0 | 4 | 0 | 1 | 5 |
| 24 | 0 | 1 | 0 | 4 | 0 | 5 |
| 25 | 1 | 0 | 0 | 0 | 4 | 5 |
| | | | | | | 625 |

| | 1 | 2 | 3 | 4 | 5 | 6 | 7 | 8 | 9 | 10 | 11 | 12 | 13 | 14 | 15 | 16 | 17 | 18 | 19 | 20 | 21 | 22 | 23 | 24 |
|---|---|---|---|---|---|---|---|---|---|---|---|---|---|---|---|---|---|---|---|---|---|---|---|---|
| 2 | 4 | | | | | | | | | | | | | | | | | | | | | | | |
| 3 | 4 | 4 | | | | | | | | | | | | | | | | | | | | | | |
| 4 | 4 | 4 | 4 | | | | | | | | | | | | | | | | | | | | | |
| 5 | 4 | 4 | 4 | 4 | | | | | | | | | | | | | | | | | | | | |
| 6 | 4 | 6 | 4 | 4 | 4 | | | | | | | | | | | | | | | | | | | |
| 7 | 4 | 4 | 4 | 6 | 4 | 4 | | | | | | | | | | | | | | | | | | |
| 8 | 4 | 4 | 4 | 4 | 6 | 8 | 8 | | | | | | | | | | | | | | | | | |
| 9 | 4 | 4 | 6 | 4 | 4 | 4 | 8 | 4 | | | | | | | | | | | | | | | | |
| 10 | 6 | 4 | 4 | 4 | 4 | 8 | 4 | 4 | 8 | | | | | | | | | | | | | | | |
| 11 | 4 | 4 | 6 | 6 | 4 | 4 | 4 | 8 | 6 | 6 | | | | | | | | | | | | | | |
| 12 | 6 | 4 | 6 | 4 | 4 | 6 | 8 | 4 | 4 | 6 | 6 | | | | | | | | | | | | | |
| 13 | 4 | 4 | 4 | 6 | 6 | 8 | 6 | 4 | 6 | 4 | 6 | 8 | | | | | | | | | | | | |
| 14 | 4 | 6 | 4 | 4 | 6 | 4 | 6 | 6 | 4 | 8 | 8 | 8 | 6 | | | | | | | | | | | |
| 15 | 6 | 6 | 4 | 4 | 4 | 6 | 4 | 6 | 8 | 4 | 8 | 6 | 8 | 6 | | | | | | | | | | |
| 16 | 4 | 6 | 4 | 8 | 8 | 6 | 4 | 6 | 8 | 6 | 4 | 10 | 4 | 6 | 8 | | | | | | | | | |
| 17 | 4 | 4 | 8 | 8 | 6 | 6 | 6 | 6 | 4 | 8 | 4 | 4 | 6 | 8 | 10 | 6 | | | | | | | | |
| 18 | 6 | 8 | 4 | 4 | 8 | 6 | 8 | 4 | 6 | 6 | 10 | 8 | 4 | 4 | 6 | 6 | 10 | | | | | | | |
| 19 | 8 | 8 | 6 | 4 | 4 | 4 | 6 | 8 | 6 | 6 | 8 | 6 | 10 | 4 | 4 | 10 | 10 | 6 | | | | | | |
| 20 | 8 | 4 | 8 | 6 | 4 | 8 | 6 | 6 | 6 | 4 | 6 | 4 | 8 | 10 | 4 | 10 | 6 | 10 | 6 | | | | | |
| 21 | 4 | 8 | 6 | 8 | 6 | 4 | 6 | 8 | 6 | 10 | 4 | 8 | 8 | 6 | 8 | 4 | 6 | 8 | 8 | 10 | | | | |
| 22 | 6 | 4 | 8 | 6 | 8 | 10 | 8 | 4 | 6 | 6 | 8 | 4 | 6 | 8 | 8 | 8 | 4 | 8 | 10 | 6 | 10 | | | |
| 23 | 8 | 6 | 4 | 6 | 8 | 8 | 6 | 6 | 10 | 4 | 8 | 8 | 4 | 8 | 6 | 6 | 10 | 4 | 8 | 8 | 10 | 8 | | |
| 24 | 6 | 8 | 8 | 4 | 6 | 6 | 10 | 6 | 4 | 8 | 8 | 6 | 8 | 4 | 8 | 10 | 8 | 6 | 4 | 8 | 8 | 8 | 10 | |
| 25 | 8 | 6 | 6 | 8 | 4 | 6 | 4 | 10 | 8 | 6 | 6 | 8 | 8 | 8 | 4 | 8 | 8 | 10 | 6 | 4 | 8 | 10 | 8 | 10 |